\documentclass{birkjour}
\usepackage{amssymb,latexsym, mathtools,tikz, hyperref}
\usetikzlibrary{arrows,decorations.markings}
\tikzstyle arrowstyle=[scale=1]
\tikzstyle directed=[postaction={decorate,
decoration={markings,mark=at position .65 with {\arrow[arrowstyle]{stealth}}}}]

\usepackage{microtype}
\usepackage{enumerate}
\usepackage{graphicx}
\usepackage{float}
\usepackage{placeins}
\usepackage{mdframed}
\usepackage{amssymb}
\usepackage{esint}
\usepackage{cool}
\usepackage[all,cmtip]{xy}
\usepackage{mathtools}
\usepackage{amstext} 
\usepackage{array}   
\usepackage[shortlabels]{enumitem}
\usepackage{ytableau}

\usepackage{cleveref}
\usepackage[utf8]{inputenc}
\usepackage{comment}

\newcolumntype{L}{>{$}l<{$}} 
\newtheorem{lemma}{Lemma}[section]
\newtheorem{theorem}[lemma]{Theorem}

\newtheorem{prop}[lemma]{Proposition}
\newtheorem{cor}[lemma]{Corollary}
\newtheorem{conj}[lemma]{Conjecture}

\newtheorem{question}[lemma]{Question}

\theoremstyle{remark}
\newtheorem{remark}[lemma]{Remark}

\theoremstyle{definition}
\newtheorem{definition}[lemma]{Definition} 
\newtheorem{setup}[lemma]{Setup}
\newtheorem{notation}[lemma]{Notation} 
 
\newtheorem{example}[lemma]{Example}

\newcommand{\inn}{\operatorname{in}}

\renewcommand{\L}{\mathcal{L}}
\newcommand{\cat}[1]{\mathcal{#1}}
\newcommand{\clique}{\Delta^{\textrm{clique}}}
\newcommand{\htt}{\operatorname{ht}}

\newcommand{\lcm}{\mathrm{lcm}}

\newcommand{\tor}{\operatorname{Tor}}

\newcommand{\pd}{\operatorname{pd}}

\renewcommand{\aa}{\mathbf a}
\newcommand{\bb}{\mathbf b}
\newcommand{\cc}{\mathbf c}
\newcommand{\dd}{\mathbf d}

\renewcommand{\O}{\mathcal{O}}

\begin{document}

\title[Determinantal Facet Ideals for Smaller Minors]{Determinantal Facet Ideals for Smaller Minors}

\author{Ayah Almousa}
\address{University of Minnesota - Twin Cities}
\email{almou007@umn.edu}
\urladdr{\url{http://umn.edu/~almou007}}

\author{Keller VandeBogert}
\address{University of Notre Dame}
\email{kvandebo@nd.edu}
\urladdr{\url{https://sites.google.com/view/kellervandebogert/}}


\urladdr{\url{http://people.math.sc.edu/kellerlv}}

\date{\today}
\begin{abstract}
A determinantal facet ideal (DFI) is generated by a subset of the maximal minors of a generic $n\times m$ matrix where $n\leq m$ indexed by the facets of a simplicial complex $\Delta$. We consider the more general notion of an $r$-DFI, which is generated by a subset of $r$-minors of a generic matrix indexed by the facets of $\Delta$ for some $1\leq r\leq n$. We define and study so-called lcm-closed and unit interval $r$-DFIs, and show that the minors parametrized by the facets of $\Delta$ form a reduced Gr\"obner basis with respect to \emph{any} term order for an lcm-closed $r$-DFI. We also see that being lcm-closed generalizes conditions previously introduced in the literature, and conjecture that in the case $r=n$, lcm-closedness is necessary for being a Gr\"obner basis. We also give conditions on the maximal cliques of $\Delta$ ensuring that lcm-closed and unit interval $r$-DFIs are Cohen-Macaulay. Finally, we conclude with a variant of a conjecture of Ene, Herzog, and Hibi on the Betti numbers of certain types of $r$-DFIs, and provide a proof of this conjecture for Cohen-Macaulay unit interval DFIs.
\end{abstract}
\maketitle

\section{Introduction}

Let $M = (x_{ij})_{1\leq i\leq n, 1\leq j\leq m}$ be an $n \times m$ matrix of indeterminates where $n\leq m$, and let $S = k[M]$ be a polynomial ring over an arbitrary field $k$ with variables in $M$.
The study of the ideal generated by all minors of a given size of $M$ has a long history, and such ideals are well understood (see, for instance, \cite{bruns2006determinantal}). In a similar vein, one can instead consider the ideal generated by \emph{some} of the minors of a given size of $M$; these are known as \emph{determinantal facet ideals} (DFIs) and were introduced by Ene, Herzog, Hibi, and Mohammadi in \cite{ene2011determinantal}. The study of DFIs turns out to be much more subtle and has seen comparably less attention, even though such ideals arise naturally in algebraic statistics (see \cite{DES98} and \cite{herzog2010binomial}). In \cite{herzog2015linear}, the linear strand of DFIs is constructed in terms of a \emph{generalized} Eagon-Northcott complex. In particular, the linear Betti numbers of such ideals may be computed in terms of the $f$-vector of an associated simplicial complex. Likewise, in \cite{vdb2020}, explicit Betti numbers of certain classes of DFIs are computed in all degrees; for arbitrary DFIs, higher degree Betti numbers have proven to be very nontrivial to compute.

DFIs for the case $n=2$ were originally introduced as binomial edge ideals independently by Ohtani \cite{ohtani2011} and Herzog, et. al. \cite{herzog2010binomial}; this generalized work of Diaconis, Eisenbud, and Sturmfels in \cite{DES98}. To study binomial edge ideals, one can associate each column of $M$ with a vertex of a graph $G$, and one can associate a minor of $M$ involving two columns $i$ and $j$ with an edge $\{i,j\}$ in the graph. For example, the ideal generated by all maximal minors of a $2\times m$ matrix corresponds to a complete graph on $m$ vertices. The relationship between homological invariants of ideals generated by some maximal minors of $M$ and combinatorial invariants of the associated graph $G$ has been widely studied; see the survey paper \cite{madani16-BEIsurvey} for a compilation of such results. DFIs naturally extend this idea by instead associating a pure simplicial complex $\Delta$ on $m$ vertices to the ideal $J_\Delta$, where each $(n-1)$-dimensional facet of $\Delta$ corresponds to a maximal minor in the set of generators of $J_\Delta$. Mohammadi and Rauh further generalized this notion to that of a determinantal hypergraph ideal, which associates a minor to each hyperedge of a graph, allowing for an ideal that is generated by minors of different sizes.

A particularly interesting class of DFIs is that for which the standard minimal generating set forms a Gr\"obner basis. It is well-known that the set of maximal minors of a generic matrix is a Gr\"obner basis for the ideal generated by them with respect to any total monomial order \cite{SturmfelsZelevinsky93, bernsteinZelevinsky93}. In the case of binomial edge ideals, there is a known necessary and sufficient condition on a graph $G$ for the generators of $J_G$ indexed by $G$ to form a reduced Gr\"obner basis with respect to a diagonal term order \cite[Theorem 1.1]{herzog2010binomial}. For a DFI where the maximal cliques of $\Delta$ overlap by $n-1$ or fewer vertices, a necessary and sufficient condition for $\Delta$ to index the generators of a reduced Gr\"obner basis for $J_\Delta$ with respect to a diagonal term order is also known \cite[Theorem 1.1]{ene2011determinantal}. 

In this paper, we introduce the notion of an $r$-DFI (Definition \ref{def: DFI}), which is an ideal generated by a subset of $r\times r$ minors in an $n\times m$ matrix indexed by the facets of some $(r-1)$-dimensional simplicial complex $\Delta$. This is a natural extension of Ene, Herzog, Hibi, and Mohammadi's DFIs in \cite{ene2011determinantal}, but is not as general as Mohammadi and Rauh's determinantal hypergraph ideals in \cite{mohammadi2018prime}.
We consider two classes of $r$-DFIs which we call \emph{lcm-closed} and \emph{unit interval} DFIs (see Definition \ref{def: lcmClosed}). In the case that $r=n$, lcm-closed DFIs are a direct generalization of some important classes of determinantal ideals in the literature, including closed binomial edge ideals and closed DFIs, but is stated for arbitrary term orders.
We apply a result of Conca (see Proposition \ref{prop:concaProp}) to show that the minimal generating set parametrized by the facets of $\Delta$ for any lcm-closed $r$-DFI forms a reduced Gr\"obner basis with respect to any diagonal term order. We then observe that in the cases that lcm-closed DFIs generalize, these conditions are also necessary for the minimal generating set to form a reduced Gr\"obner basis. This leads us to pose Question \ref{q: lcmClosedNecessary}, which asks whether or not being lcm-closed is a necessary condition for being a Gr\"obner basis when $r=n$. In the case $n=2$, we prove that Question \ref{q: lcmClosedNecessary} is true.

In Section \ref{sec:CMness}, we apply the results of Section \ref{sec: lcmClosed} to deduce certain cases for which lcm-closed and unit interval DFIs must be Cohen-Macaulay.
We start with a result that allows one to deduce when the sum of two Cohen-Macaulay ideals remains Cohen-Macaulay based on knowledge of their initial ideals (see Proposition \ref{prop:CMnessCondition}). This leads to Corollary \ref{cor:CMnessDFI}, which proves that if the maximal cliques of $\Delta$ have sufficiently small pairwise intersections, then the associated $r$-DFI must be Cohen-Macaulay.
Lastly, we conclude with a variant of a conjecture by Ene, Herzog, and Hibi (see Conjecture \ref{conj:equalBnos}) and give a proof of the conjecture in the case that $J_\Delta$ is a Cohen-Macaulay unit interval $r$-DFI.

\section{Lcm-closed and Unit Interval Determinantal Facet Ideals}\label{sec: lcmClosed}
In this section, we generalize the idea of a determinantal facet ideal to that of an $r$-determinantal facet ideal. We recall examples in the literature, some of which are specializations of Definition \ref{def: DFI}, and introduce a sufficient condition for the standard minimal generating set of these so called $r$-determinantal facet ideals to be a reduced Gr\"obner basis (see Definition \ref{def: lcmClosed}). This condition generalizes conditions considered in the binomial edge ideal case and for that of closed determinantal ideals (see \cite{ene2011determinantal}). Let us set the stage with some notation and definitions:

\begin{notation}\label{not: intro}
Fix $r$ to be a positive integer. Let $M = (x_{ij})_{1\leq i\leq n, 1\leq j\leq m}$ be an $n \times m$ matrix of indeterminates where $n\leq m$, and let $S = k[M]$ be a polynomial ring over an arbitrary field $k$ with variables in $M$.
For indices $\aa = \{a_1,\ldots, a_r\}$ and $\bb = \{b_1, \ldots, b_r\}$ such that $1\leq a_1 < \cdots < a_r\leq n$ and $1\leq b_1 < \cdots < b_r\leq m$, set
$$
[\aa\vert\bb] = [a_1,\ldots, a_r\vert b_1,\ldots, b_r] = \det \left( \begin{array}{ccc}
    x_{a_1,b_1} & \cdots & x_{a_1,b_r} \\
    \vdots & \ddots & \vdots\\
    x_{a_r,b_1} & \cdots & x_{a_r,b_r}\\
\end{array} \right)
$$
where $[\aa\vert\bb]=0$ if $r>n$. When $r=n$, use the simplified notation $[\aa]$ = $[1,\ldots, n\vert \aa]$. The ideal generated by the $r$-minors of $M$ is denoted $I_r(M)$.
\end{notation}

\begin{definition}\label{def:cliques}
Let $\Delta$ be a pure $(r-1)$-dimensional simplicial complex on vertex set $[m]$. For an integer $i$, the $i$-th skeleton $\Delta^{(i)}$ of $\Delta$ is the subcomplex of $\Delta$ whose faces are those faces of $\Delta$ with dimension at most $i$. Let $\cat{S}$ denote the set of simplices $\Gamma$ with vertices in $[m]$ with $\dim (\Gamma) \geq r-1$ and $\Gamma^{(r-1)} \subset \Delta$. 

Let $\Gamma_1 , \dotsc , \Gamma_c$ be maximal elements in $\cat{S}$ with respect to inclusion, and let $\Delta_i := \Gamma^{(r-1)}_i$. Each $\Gamma_i$ is called a \emph{maximal clique}, and any induced subcomplex of $\Gamma_i$ is a \emph{clique}. The simplicial complex $\clique$ whose facets are the maximal cliques of $\Delta$ is called the \emph{clique complex} associated to $\Delta$. The decomposition $\Delta = \Delta_1 \cup \cdots \cup \Delta_c$ is called the \emph{maximal clique decomposition} of $\Delta$.
\end{definition}

\begin{definition}\label{def: DFI}
Adopt Notation \ref{not: intro}, and let $\Delta$ be an $(r-1)$-dimensional simplicial complex on the vertex set $[m]$.
The \emph{$r$-determinantal facet ideal} (or \emph{$r$-DFI}) $J_\Delta\subseteq S$ associated to $\Delta$ is the ideal generated by determinants of the form $[\aa | \bb]$ where $\bb$ supports an $(r-1)$-face of $\Delta$; that is, the columns of $[\aa | \bb]$ correspond to the vertices of some $(r-1)$-face of $ \Delta$.
\end{definition}

\begin{notation}\label{not: rDFI}
Let $\Delta$ be a pure $(r-1)$-dimensional simplicial complex on the vertex set $[m]$ with maximal clique decomposition $\Delta = \Delta_1 \cup \cdots \cup \Delta_c$. The notation $J_{\Delta_i}$ will be used to denote the $r$-DFI associated to the simplicial complex $\Delta_i$.
\end{notation}

\begin{remark}
Definition \ref{def: DFI} naturally generalizes the notion of a DFI introduced by Ene, Herzog, Hibi, and Mohammadi in \cite{ene2011determinantal}, who considered only the case when $r=n$. DFIs are, in turn, a generalization of binomial edge ideals introduced in \cite{herzog2010binomial}, which coincides with the case when $r=n=2$ and $\Delta$ is a graph $G$. However, $r$-DFIs are not as general as Mohammadi and Rauh's notion of a \textit{determinantal hypergraph ideal} introduced in \cite{mohammadi2018prime}, which allows for the ideal to be generated by minors of different sizes in $M$.
\end{remark}

\begin{remark} Let $J_\Delta$ denote any $r$-DFI. The simplicial complex $\Delta$ associated to an $r$-DFI serves as a combinatorial tool to parametrize the column sets appearing on minimal generators of $J_\Delta$. Maximal cliques in the clique decomposition of $\Delta = \bigcup_{i=1}^c \Delta_i$ correspond to the largest submatrices $M_i$ of $M$ such that the ideal generated by all $r$-minors of $M_i$ is contained in $J_\Delta$.
\end{remark}

\begin{example}\label{ex: DFI} 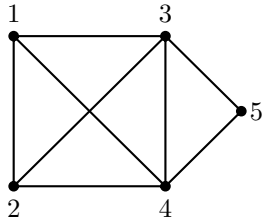
\begin{figure}[ht]
\begin{tikzpicture}[scale=0.4]
\filldraw[thick] 
(0,0) circle (4pt) --
(5,0) circle (4pt)  -- 
(5,5) circle (4pt) --
(0,5) circle (4pt) --
(0,0);
\filldraw[thick] (5,5)--(7.5, 2.5) circle (4pt) -- (5,0);
\draw[thick] (0,0)--(5,5);
\draw[thick] (0,5)--(5,0);
\node at (-0,-0.75) {$2$};
\node at (5,-0.75) {$4$};
\node at (5,5.75) {$3$};
\node at (-0,5.75) {$1$};
\node at (8, 2.5) {$5$};
\end{tikzpicture}
\caption{A graph with two maximal cliques given by vertex sets $\{1,2,3,4\}$ and $\{3,4,5\}$.}\label{fig: BEI2cliques}
\end{figure} Let $G$ be the graph in Figure \ref{fig: BEI2cliques}. If $r = n = 2$, then $J_G$ corresponds to a subideal of the ideal of maximal minors of a $2\times 5$ matrix with generators $[a_1 a_2]$ indexed by the $(a_1,a_2)$ in the edge set of $G$.

If, instead, $n=3$, then $J_G$ corresponds to a subideal of the ideal of $I_2(M)$ where $M$ is a $3\times 5$ matrix. Now the generators of $J_G$ are of the form $[a_1 ,a_2 | b_1, b_2]$ where $\{a_1,a_2\}\subset [3]$ and $\{b_1<b_2\}$ is an edge of $G$.
\end{example}

\begin{notation}\label{not: leadTerm} Let $<$ be a total monomial order in a polynomial ring $S$ over a field $k$. If $f$ is a polynomial in $S$, denote by $\inn_<(f)$ the leading term of $f$ with respect to $<$. If $I\subseteq S$ is an ideal, denote by $\inn_<(I)$ the initial ideal of $I$ with respect to $<$. Frequently, when the term order $<$ is clear, it will be dropped and leading terms and initial ideals will simply be denoted by $\inn(f)$ and $\inn(I)$, respectively.
\end{notation}



The following Proposition, originally due to Conca, will turn out to be very useful in the proofs of Theorems \ref{thm: lcmGB} and \ref{thm:intervalsGB}.

\begin{prop}[{\cite[Lemma 1.3]{conca1996gorenstein}}]\label{prop:concaProp} Let $S$ be a polynomial ring over a field $k$, and let $<$ be a term order. Let $I$ and $J$ be homogeneous ideals of $S$. Then
\begin{enumerate}[label = \alph*)]
    \item $\inn(I)+\inn(J)\subseteq \inn(I+J)$ and $\inn(I\cap J)\subseteq \inn(I)\cap \inn(J)$,
    \item $\inn(I)+\inn(J) = \inn(I+J)$ if and only if $\inn(I\cap J) = \inn(I)\cap \inn(J)$,
    \item let $F$ be a Gr\"obner basis of $I$ and let $G$ be a Gr\"obner basis of $J$. Then $F\cup G$ is a Gr\"obner basis of $I+J$ if and only if, for all $f\in F$ and $g\in G$, there exists $h\in I\cap J$ such that $\inn(h) = \lcm(\inn(f),\inn(g))$.
\end{enumerate}
\end{prop}


\begin{definition}\label{def: lcmClosed} Adopt Notation \ref{not: intro} and Notation \ref{not: leadTerm}. Let $<$ be any term order. Let $\Delta$ be a pure $(r-1)$-dimensional simplicial complex on $m$ vertices with maximal clique decomposition $\Delta = \bigcup_{i=1}^c \Delta_i$. The $r$-DFI $J_\Delta$ is \emph{lcm-closed} if the following condition holds:
\begin{enumerate}[(*)]
    \item For all $[\aa | \bb] \in J_{\Delta_i}$, $[\aa' | \bb'] \in J_{\Delta_j}$ with non-coprime lead terms and $[\aa | \bb], \ [\aa' | \bb'] \notin J_{\Delta_i \cap \Delta_j}$, there exists $[\cc | \dd] \in J_{\Delta_i \cap \Delta_j}$ such that $\inn_<([\cc | \dd])$ divides $\lcm\big( \inn_<([\aa | \bb] ) , \inn_<( [\aa' | \bb']) \big)$. 
\end{enumerate}
The $r$-DFI $J_\Delta$ is a \emph{unit interval DFI} if each $\Delta_i$ may be written as an interval $[a_i , b_i] = \{ a_i , a_i+1 , \dots , b_i-1 , b_i \}$ for integers $a_i < b_i$.
\end{definition}

\begin{example} Let $\Delta$ be a $2$-dimensional simplicial complex on $5$ vertices with clique decomposition $\Delta_1 = \{1,2,3,4\}$ and $\Delta_2 = \{2,3,4,5\}$, and let $r = n = 3$. Let $<$ be any diagonal term order. For any $2$-dimensional faces $\aa\in\Delta_1$ and $\aa' \in \Delta_2$ such that $\aa,\aa'\notin \Delta_1\cap\Delta_2$, observe that $\inn_<[\aa]$ and $\inn_<[\aa']$ are coprime except for the case where $\aa = [1,3,4]$ and $\aa' = [2,3,5]$. In this case,
\begin{align*}
\lcm(\inn_<([1,3,4]), \inn_<([2,3,5]) &= \lcm(x_{11} x_{23} x_{34}, x_{12} x_{23} x_{35}) \\
&= x_{11} x_{12} x_{23} x_{34} x_{35}
\end{align*}
which is divisible by $x_{12} x_{23} x_{34} = \inn_<([2,3,4])$, and $\inn_<([2,3,4]) = J_{\Delta_1\cap\Delta_2}$. Therefore, $J_\Delta$ is lcm-closed and, in particular, is a unit interval $3$-DFI.
\end{example}

The following classes of $r$-DFIs are lcm-closed. Each of these classes forms a Gr\"obner basis when $<$ is the standard diagonal term order.

\begin{example}[Closed binomial edge ideals]\label{ex: closedBEI} When $r=n=2$, $\Delta$ may be associated with a graph $G$ and $J_G$ is called a binomial edge ideal. Denote by $E(G)$ the edge set of $G$. A graph $G$ (or its respective binomial edge ideal $J_G$) is \emph{closed} with respect to a given labeling if, for all distinct pairs of edges $\{i,j\}$ and $\{k,\ell\}$ with $i<j$, $k<\ell$, one has $\{j,\ell\}\in E(G)$ if $i=k$, and $\{i,k\}\in E(G)$ if $j=\ell$. This condition is equivalent to $J_G$ having a quadratic Gröbner basis \cite[Theorem 1.1]{herzog2010binomial}.

If two generators $\inn_<[\aa]$ and $\inn_<[\bb]$ of $J_G$ are not relatively prime, then $J_G$ satisfies (*) if and only if they correspond to edges $\aa$ and $\bb$ in the same clique of $G$. This is equivalent to the closed condition for a binomial edge ideal, so $J_G$ is lcm-closed if and only if it is closed.
\end{example}

\begin{example}[Closed DFIs]\label{ex: closedDFI} Let $r=n$ and assume that a pure $n$-dimensional simplicial complex $\Delta$ satisfies the following condition:
\begin{itemize}
\item[$\blacklozenge$] \text{No two maximal cliques of $\Delta$ share more than $n-1$ vertices.}
\end{itemize}
A DFI satisfying $\blacklozenge$ is said to be \textit{closed} if for all $i\neq j$ and all $\{a_1 < \dots < a_n\}\in \Delta_i$ and $\{b_1 < \dots < b_n\}\in \Delta_j$, the monomials $\inn_<[a_1,\dots, a_n]$ and $\inn_<[b_1,\dots, b_n]$ are relatively prime. A DFI $J_\Delta$ satisfying $\blacklozenge$ is closed if and only if  generating the generating set of $J_\Delta$ indexed by the facets of $\Delta$ forming a Gröbner basis with respect to $<$ \cite[Theorem 1.1]{ene2011determinantal}.

If $\Delta$ satisfies $\blacklozenge$, then there are no generators of $J_\Delta$ contained in the intersection of any two maximal cliques, so $J_\Delta$ satisfies (*) if and only if it is closed.
\end{example}

\begin{setup}\label{set: lcmClosed} Let $\Delta$ be a pure $(r-1)$-dimensional simplicial complex on the vertex set $[m]$ admitting maximal clique decomposition $\Delta = \bigcup_{i=1}^c \Delta_i$. Let $S = k[x_{ij} \mid 1 \leq i \leq n, \ 1 \leq j \leq m]$ be a polynomial ring over an arbitrary field $k$.
\end{setup}

\begin{theorem}\label{thm: lcmGB}
Adopt notation and hypotheses as in Setup \ref{set: lcmClosed}. If the associated $r$-DFI $J_\Delta$ is lcm-closed, then the generators of $J_\Delta$ indexed by the facets of $\Delta$ form a reduced Gr\"obner basis with respect to any term order $<$.
\end{theorem}

\begin{proof}
Assume $ f = [\aa\mid\bb]$ is a minimal generator of $J_{\Delta_i}$ and  $ g = [\aa'\mid\bb']$ is a minimal generator of $J_{\Delta_j}$ such that $f,g\notin J_{\Delta_i}\cap J_{\Delta_j}$. By Proposition \ref{prop:concaProp}, it suffices to check that there exists some $h\in J_{\Delta_i}\cap J_{\Delta_j}$ such that $\inn(h) = \lcm(\inn(f),\inn(g))$. If $\inn(f)$ and $\inn(g)$ are coprime, then 
\begin{align*}
\lcm(\inn(f),\inn(g)) &= \inn(f)\cdot \inn(g) \\
&= \inn(f\cdot g)\in J_{\Delta_i}\cdot J_{\Delta_j}\subseteq J_{\Delta_i}\cap J_{\Delta_j}.
\end{align*}
Suppose that $\inn(f)$ and $\inn(g)$ are not coprime. By the definition of lcm-closed, there exists some $h\in J_{\Delta_i\cap \Delta_j}\subseteq J_{\Delta_i}\cap J_{\Delta_j}$ such that $\inn(h)$ divides $\lcm(\inn(f),\inn(g))$. Any multiple of $h$ is also contained in both $J_{\Delta_i}$ and $J_{\Delta_j}$, so $\frac{\lcm(\inn(f),\inn(g))}{\inn(h)}\cdot h \in J_{\Delta_i}\cap J_{\Delta_j}$ and has initial term equal to $\lcm(\inn(f),\inn(g))$. 
\end{proof}

\begin{example}\label{ex:intervalsexample}
In the case $r=n$ and $<$ is a diagonal term order, it is clear that unit interval DFIs are lcm-closed. For $r<n$, these notions are distinct. Let $n=3$ and consider for instance the $2$-DFI $J_\Delta$ associated to the simplicial complex $\Delta$ with maximal clique decomposition $[1,2] \cup [2,3]$. Notice that $J_\Delta$ is \emph{not} lcm-closed with respect to a diagonal term order because the generators $[12 | 12]$ and $[23 | 23]$ do not have coprime lead terms. However, $J_\Delta$ is a unit interval $r$-DFI and its natural minimal generating set forms a reduced Gr\"obner basis with respect to any diagonal term order. 
In particular, one finds that the determinant of
$$\begin{pmatrix}
      {x}_{11}&{x}_{12}&{x}_{13}\\
      {x}_{21}&{x}_{22}&{x}_{23}\\
      {x}_{31}&{x}_{32}&{x}_{33}\end{pmatrix}$$
is an element of $J_{[1,2]} \cap J_{[2,3]}$ whose lead term divides the lcm of the lead terms of $[12|12]$ and $[23|23]$.
\end{example}

The following theorem can be proved using rather computationally intense methods, but follows much more easily from work of Seccia (see \cite{seccia2021knutson}). 

\begin{theorem}\label{thm:intervalsGB}
Adopt notation and hypotheses as in Setup \ref{set: lcmClosed}. If $J_\Delta$ is a unit interval $r$-DFI, then the generators of $J_\Delta$ indexed by the facets of $\Delta$ form a reduced Gr\"obner basis with respect to any diagonal term order $<$.
\end{theorem}

\begin{proof} See Corollary 2.4 of \cite{seccia2021knutson}.
\end{proof}



As shown in Examples \ref{ex: closedBEI} and \ref{ex: closedDFI}, lcm-closed $r$-DFIs coincide with many known cases of determinantal ideals which form a Gr\"obner basis in the literature. Moreover, in all of the cases considered previously, the property of being lcm-closed is necessary for the standard minimal generating set to form a Gr\"obner basis. This leads us to ask Question \ref{q: lcmClosedNecessary}.

\begin{question}\label{q: lcmClosedNecessary} Let $\Delta$ be a pure $(r-1)$-dimensional simplicial complex. Is the property of being lcm-closed not only sufficient, but also necessary, for the generators of $J_\Delta$ indexed by the facets of $\Delta$ to form a reduced Gr\"obner basis with respect to any term order $<$?
\end{question}

As evidence of a positive answer to Question \ref{q: lcmClosedNecessary}, we have the following proposition:

\begin{prop}
Question \ref{q: lcmClosedNecessary} is true for binomial edge ideals.
\end{prop}

\begin{proof}
Let $<_1$ be any term order on $S$. By Proposition $1.11$ of \cite{SturmfelsZelevinsky93}, there exists a permutation $\sigma : [m] \to [m]$ such that $\inn_{<_1} ([a_1,a_2]) = \inn_{<_2} ([\sigma(a_1) , \sigma(a_2) ])$, where $\{ a_1 , a_2 \} \in \Delta$ and $<_2$ denotes the standard diagonal term order. Observe that $\sigma$ induces an automorphism of $k$-algebras $\phi : S \to S$ by acting on the second indices of the variables. By Theorem $2.2$ of \cite{ene2011cohen}, the ideal $\phi (J_\Delta) = J_{\sigma(\Delta)}$ is lcm-closed with respect to $<_2$. Since the definition of lcm-closed is an algebraic condition, it follows that $J_\Delta$ is lcm-closed with respect to $<_1$.
\end{proof}

\section{Cohen-Macaulayness for Certain Classes of $r$-DFIs}\label{sec:CMness}

Recall that an ideal $I$ in a standard graded polynomial ring $S$ over a field $k$ is \emph{Cohen-Macaulay} if $\htt (I) = \pd_S (S/I)$. In this section, we investigate Cohen-Macaulayness for certain classes of $r$-DFIs. In particular, we give sufficient conditions for the Cohen-Macaulayness of both lcm-closed and unit interval $r$-DFIs in Corollary \ref{cor:CMnessDFI}. This allows us to prove a variant of a conjecture of Ene, Herzog, and Hibi for Cohen-Macaulay unit interval DFIs (Corollary \ref{cor:intervalBnos}), and we pose a similar conjecture for $r$-DFIs in which each clique has precisely $r$ vertices.

The following Proposition is likely well-known, but does not seem to appear explicitly in the literature. We give a complete statement and proof for reference and convenience; it turns out to be a surprisingly effective method for deducing Cohen-Macaulayness of sums of ideals.

\begin{prop}\label{prop:CMnessCondition}
Let $S$ be a standard graded polynomial ring over a field $k$ and let $<$ be any term order. Assume that $I$ and $J$ are Cohen-Macaulay ideals with the property that $$\inn_< (I \cap J) = \inn_< (I) \cap \inn_< (J) = \inn_< (I) \cdot \inn_< (J).$$
Then $I+J$ is a Cohen-Macaulay ideal.
\end{prop}

\begin{proof}
It is a standard fact that $\htt (I) = \htt ( \inn_< (I))$ for any ideal $I$. One has the following string of equalities:
\begingroup\allowdisplaybreaks
\begin{align*}
    \htt (I+ J ) &= \htt ( \inn_< (I + J) ) \\
    &= \htt (\inn_< (I) + \inn_< (J)) \qquad \textrm{(by Proposition \ref{prop:concaProp})} \\
    &= \htt(\inn_< (I)) + \htt (\inn_< (J)) \quad \textrm{(since} \ \inn_< (I) \cap \inn_< (J) = \inn_< (I) \cdot \inn_< (J) \textrm{)} \\
    &= \htt (I) + \htt (J). 
\end{align*}
\endgroup
Likewise, there is a string of implications:
\begingroup\allowdisplaybreaks
\begin{align*}
    &\inn_< (I) \cap \inn_< (J) = \inn_< (I) \cdot \inn_< (J) \\ \implies& I \cap J = IJ \\
    \iff& \tor_1^S (S / I , S/ J ) = 0 \\
    \iff& \tor_i^S (S/I , S/J) = 0 \ \textrm{for all} \ i>0 \quad \textrm{(by \cite[Corollary 2.5]{celikbas2015vanishing})} \\
    \implies& \pd_S (S/(I+J)) = \pd_S (S/I) + \pd_S (S/J). 
\end{align*}
\endgroup
Since $I$ and $J$ were assumed to be Cohen-Macaulay, the result follows.
\end{proof}

Recall that for a generic $n \times m$ matrix $M$, the ideal of $r$-minors $I_r (M)$ is Cohen-Macaulay for any $1 \leq r \leq n$ (see \cite{bruns2006determinantal}). This implies that for each maximal clique $\Delta_i$ appearing in the clique decomposition of $\Delta$, $J_{\Delta_i}$ is a Cohen-Macaulay ideal; in particular, one uses this to deduce the following:

\begin{cor}\label{cor:CMnessDFI}
Adopt notation and hypotheses as in Setup \ref{set: lcmClosed}. Then:
\begin{enumerate}
    \item If $J_\Delta$ is an lcm-closed $r$-DFI with $<$ any term order and $|V(\Delta_i) \cap V(\Delta_j)| \leq r-1$, then $J_\Delta$ is Cohen-Macaulay.
    \item If $J_\Delta$ is a unit interval $r$-DFI with $<$ any diagonal term order and $|V(\Delta_i) \cap V(\Delta_j)| \leq \max \{ 0 , 2r-n-1 \}$, then $J_\Delta$ is Cohen-Macaulay.
\end{enumerate}
\end{cor}

\begin{proof}
Both of the conditions on the intersection size of each clique in $(1)$ and $(2)$ ensure that minimal generators coming from any two distinct maximal cliques have coprime lead terms. The result then follows immediately upon combining Theorems \ref{thm: lcmGB} and \ref{thm:intervalsGB} with Proposition \ref{prop:CMnessCondition}.
\end{proof}

The following Corollary yields a proof of a variant of a conjecture by Herzog and Hibi (see \cite{ene2011cohen}) for a certain class of unit interval DFIs. It states that the graded Betti numbers of $J_\Delta$ and its initial ideal with respect to $<$ are always equal. In general, the Betti numbers of the initial ideal are only an upper bound; it is \emph{very} rare to have equality everywhere.

\begin{cor}\label{cor:intervalBnos}
Adopt notation and hypotheses as in Setup \ref{set: lcmClosed} with $<$ any diagonal term order and assume that $J_\Delta$ is a unit interval $r$-DFI with $|V(\Delta_i) \cap V(\Delta_j)| \leq \max \{ 0 , 2r-n-1 \}$. If $|V (\Delta_i )| = r$ for each $i=1, \dots , c$, then
$$\beta_{ij} (S/ J_\Delta) = \beta_{ij} (S / \inn_< (J_\Delta)) \ \textrm{for all} \ i,j.$$
\end{cor}

\begin{proof}
Let $F^{\Delta_i}_\bullet$ denote the minimal free resolution of each $S/ J_{\Delta_i}$. The condition $|V(\Delta_i) \cap V(\Delta_j)| \leq \max \{ 0 , 2r-n-1 \}$ implies that minimal generators coming from any two distinct maximal cliques have coprime lead terms, whence the minimal free resolution of $S/J_\Delta$ may be obtained as the tensor product complex $F^{\Delta_1}_\bullet \otimes \cdots \otimes F^{\Delta_c}_\bullet$. Since $|V(\Delta_i)| = r$, one has $\beta_{jk} (S/J_{\Delta_i} ) = \beta_{jk} (S / \inn_< J_{\Delta_i})$ for all $j,k$, where $i = 1 , \dots , c$ (see, for instance, \cite[Theorem 1.4]{boocher2012}). Combining the previous two sentences yields the result.
\end{proof}

Corollary \ref{cor:intervalBnos} combined with copious amounts of computational evidence suggests that the following conjecture holds:

\begin{conj}\label{conj:equalBnos}
Adopt notation and hypotheses as in Setup \ref{set: lcmClosed} and assume that the standard minimal generating set of $J_\Delta$ forms a reduced Gr\"obner basis. If $|V(\Delta_i) | = r$ for all $i = 1 , \dots , c$, then
$$\beta_{ij} (S/J_\Delta) = \beta_{ij} (S/ \inn_< (J_\Delta)) \ \textrm{for all} \ i,j.$$
\end{conj}

\begin{example}\label{ex:nonCMequalBnos}
Let $n=4$ and consider the $3$-DFI associated to the simplicial complex with maximal clique decomposition $[1,3] \cup [2,4]$. Let $<$ be a diagonal term order. It can be shown using Macaulay2 \cite{M2} that $J_\Delta$ is \emph{not} Cohen-Macaulay and that $S/ J_\Delta$ and $S/ \inn_< (J_\Delta)$ both have Betti table
$$\begin{matrix}
       &0&1&2&3&4\\\text{total:}&1&8&17&16&6\\\text{0:}&1&\text{.}&\text{.}&
       \text{.}&\text{.}\\\text{1:}&\text{.}&\text{.}&\text{.}&\text{.}&\text
       {.}\\\text{2:}&\text{.}&8&7&\text{.}&\text{.}\\\text{3:}&\text{.}&\text
       {.}&\text{.}&\text{.}&\text{.}\\\text{4:}&\text{.}&\text{.}&10&16&6.\\
       \end{matrix}$$
\end{example}
\section*{Acknowledgments} We would like to thank the anonymous referee for their careful reading and helpful comments, especially for pointing out an error in our original proof of Theorem \ref{thm:intervalsGB}. We would also like to thank Matteo Varbaro for informing us of the paper \cite{seccia2021knutson} by Seccia. The first author was partially supported by the NSF GRFP under Grant No. DGE-1650441.

\bibliographystyle{amsplain}
\bibliography{biblio}
\addcontentsline{toc}{section}{Bibliography}

\end{document}